# Functional approach for excess mass estimation in the density model

## Cristina Butucea

*Laboratoire Paul Painlevé - UMR CNRS 8524*
*Université Lille 1*
*59655 Villeneuve d'Ascq CEDEX, France*
*e-mail:* `cristina.butucea@math.univ-lille1.fr`

## Mathilde Mougeot and Karine Tribouley

*Modal'X, Université Paris X,*
*200, avenue de la République*
*92001 Nanterre Cedex, France*
*and Laboratoire de probabilités et modèles aléatoires,*
*Université Paris VI, 4, place Jussieu, Boîte courrier 188,*
*75252 Paris, France*
*e-mail:* `mathilde.mougeot@u-paris10.fr`
*e-mail:* `karine.tribouley@u-paris10.fr`

**Abstract:** We consider a multivariate density model where we estimate the excess mass of the unknown probability density $f$ at a given level $\nu > 0$ from $n$ i.i.d. observed random variables. This problem has several applications such as multimodality testing, density contour clustering, anomaly detection, classification and so on. For the first time in the literature we estimate the excess mass as an integrated functional of the unknown density $f$. We suggest an estimator and evaluate its rate of convergence, when $f$ belongs to general Besov smoothness classes, for several risk measures. A particular care is devoted to implementation and numerical study of the studied procedure. It appears that our procedure improves the plug-in estimator of the excess mass.

## 1. Introduction

Let $X_1, \ldots, X_n$ be $n$ i.i.d. observations in $\mathbb{R}^d$, $d \geq 1$ having unknown underlying distribution function $F$ with probability density $f$. We want to estimate the excess mass of this distribution, at level $\nu > 0$ which was defined by (19) as

$$\mathcal{E}(\nu) = F(C(\nu)) - \nu \cdot |C(\nu)|,$$

where $|\cdot|$ denotes the Lebesgue measure of a set and $C(\nu) = \{x \in \mathbb{R}^d : f(x) \geq \nu\}$ is the density level set (at level $\nu$) or density contour cluster (see Figure 1).

449

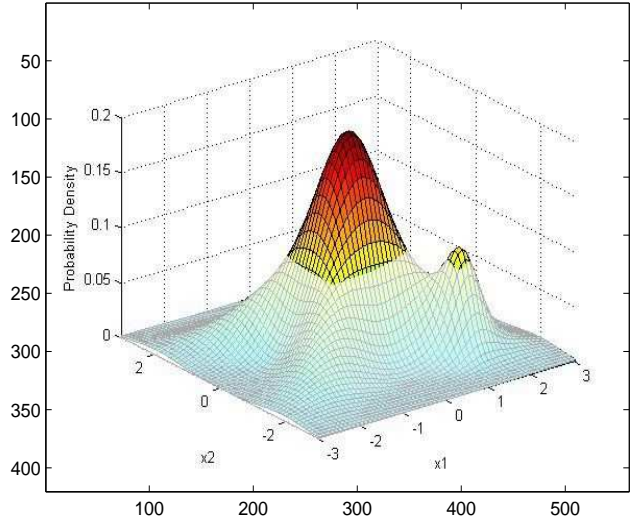

FIG 1. *Excess mass for a bivariate probability density with two local modes*

Estimating the excess mass has multiple practical applications that we mention without actually dealing with them. Most applications use differences of excess-masses at different levels $\nu$ in order to test the multimodality of a probability distribution. Hartigan and Hartigan (11) introduced the dip-excess mass and defined an estimator which allowed to test multimodality. This estimator was extensively used in the literature since, see e.g. (19), (5), (7). They insist on the fact that such a procedure separates mode estimation from its location.

Another important application of the excess mass functional is to the estimation density level sets (or density contour clustering), that is the support (the set of points) $C(\nu)$ on which the excess mass at level $\nu$ is calculated. This requires a good estimator of the excess mass as well as an optimization procedure. Polonik (21) proved consistency of such estimators of the density level set and found some rates of convergence. Tsybakov (26) gave minimax rates for estimating smooth star-shaped level sets of a density. These methods are either very difficult to implement or use assumptions which are difficult to check. They use a margin assumption quantifying the smoothness of the density $f$ around the level $\nu$ as introduced by (17). Later, (8) used a Bayesian approach and (23) revisited the plug-in estimator for this problem. They may claim for computational feasibility as well as for strong theoretical properties. On the other hand, (14) studied and implemented an estimator of the support of a density via complexity penalized excess-mass criterion.

Other applications of excess mass estimation include estimation of regression contour clusters (22), discrimination of locally stationary time series (4) and, via level set estimation, anomaly detection and classification as described by (23).

These methods generally avoid using a nonparametric estimator of $f$. Indeed, such an estimator may not be very attractive in higher dimensions $d$. We overpass this difficulty by estimating the excess mass $\mathcal{E}(\nu)$ as an integrated functional of $f$ at fixed level $\nu > 0$, that is

$$\mathcal{E}(\nu) = \int \Phi_\nu(f(t))dt \qquad \text{for} \qquad \Phi_\nu(x) = (|x| - \nu) \, 1_{|x| > \nu}. \qquad (1.1)$$

Indeed, excess mass estimation is a particular case of estimating integrated functionals of $f$ of general type: $\theta = \int \Phi(f)$, where $\Phi$ is known. The study of such functionals with $\Phi$ 4-times continuously differentiable is now completed. It was noticed since (1), (12) and many others that $\theta$ can be estimated at a parametric rate as soon as the Hölder smoothness of $f$ is larger or equal to $1/4$, but at a slower nonparametric rate otherwise. The lower bounds for the nonparametric rates case were established in (2). These rates were achieved in a minimax setup by wavelet estimation procedure in the paper by (13) and in an adaptive to the smoothness setup in the paper by (24) (with a loss with respect to the minimax rate of the usual logarithmic order). Nemirovski (20) gave asymptotically efficient estimators for 1 and 2-times continuously differentiable function $\Phi$. In our problem $\Phi$ is continuous but not differentiable (when periodized). Our approach works for any other integrated functionals with continuous but not differentiable $\Phi$.

In the particular case of a large enough level $\nu$, the excess mass problem is reduced to the estimation of the $\mathbb{L}_1$ norm. Obviously, this problem has no interest in the density model. In the regression model (15) studied the problem of estimating the $\mathbb{L}_1$ norm and in the gaussian white noise model (16) estimated the $\mathbb{L}_r$ norm for $r \geq 1$.

The excess mass estimator we construct in this paper generalizes the estimator of the $\mathbb{L}_1$ norm in (16). Their procedure actually uses a Fourier series approximation for the function $\Phi$, whose coefficients are known and depend on the level $\nu$. As $\Phi$ is applied to $f$, so are the functions of the Fourier basis. A kernel estimator of $f$ is then plugged-into the functions of the Fourier basis and they are multiplied by a factor which actually reduces the bias. This multiplicative factor is depending on the variance of the kernel estimator in the considered model. The integral of this expansion gives the final estimator of the functional.

In this paper, we consider a density (thus heteroscedastic) model. Nevertheless, the excess mass functional can be defined for the regression and gaussian white noise models as well. As we consider the density model, the multiplicative factor depends upon a variance which is proportional to the unknown density $f$ and another estimator is plugged-into this factor. Moreover, we have a multidimensional setup fitting better to most applications. For the study of rates, we replace the preliminary kernel estimator by a wavelet estimator which allows us to compute rates over more general Besov smoothness classes for the unknown probability density $f$. The whole procedure is detailed and fully explained in Section 2.

In (16), lower bounds for estimating the $\mathbb{L}_1$-norm were given, but upper and lower bounds were separated by a logarithmic factor. We show that our estimator attains the same rate as the estimator of the $\mathbb{L}_1$-norm in (16). In the gaussian white noise model, (3) improved on the lower bounds for the particular problem of excess mass estimation. Nevertheless, a gap still remains between the upper bounds we present here and their lower bounds.

The paper is organized as follows. In Section 2, we present the estimation procedure. In Section 3, we state an expansion for the upper bound of the expected errors (pointwise, $\mathbb{L}_2$ and $\mathbb{L}_\infty$) of our procedure. Next, we determine the optimal parameters of the method and give in Theorem 3.1 the rate our procedure achieves. Section 4 is devoted to the empirical study. The proofs are postponed to Section 5.

## 2. Excess mass estimation procedure

Let us describe the densities $f$ considered in the sequel. We suppose that $f$ is compactly supported with known support $\mathbb{K} = [A_1, B_1] \times \ldots \times [A_d, B_d] \subseteq \mathbb{R}^d$. We denote, for some $m^* > 0$, $\mathcal{F}(\mathbb{K}, m^*)$ the class of compactly supported probability densities $f : \mathbb{K} \to \mathbb{R}$ such that

$$\inf_{t \in \mathbb{K}} f(t) \geq m^* \quad \text{and} \quad \|f\|_\infty \leq \rho \tag{2.1}$$

is satisfied for some $\rho \in ]0, 1[$.

The estimation procedure consists of four steps. At the first step we approximate the functional $\Phi_\nu$ defined in (1.1) by its truncated Fourier series with known coefficients. Then, at the second step, we estimate the unknown function $f$. In particular, we consider here wavelet estimator but it is possible to consider kernel estimators. The third step consists in plugging-into the Fourier series the wavelet estimator of $f$ in an unbiased way. Finally, we integrate on $\mathbb{K}$ to get the estimator $\hat{\mathcal{E}}(\nu)$ of $\mathcal{E}(\nu)$. Let us describe in more details this procedure.

### 2.1. Approximation of the functional $\Phi_\nu$

We assumed in (2.1) that the density of interest $f$ is bounded by some $\rho < 1$ uniformly over the class $\mathcal{F}(\mathbb{K}, m^*)$. It allows us to define $\Phi_\nu$ as a function on $[-1, 1]$ to $[0, 1]$ and then, its approximation by Fourier series is given by

$$A_N \Phi_\nu(u) \stackrel{\text{def}}{=} c_0(\nu) + \sum_{k=1}^{N} (c_k(\nu) \cos(\pi k u) + b_k(\nu) \sin(\pi k u)),$$

where the Fourier coefficients are easily computed

$$\begin{cases} c_0(\nu) &= \langle 1, \Phi_\nu \rangle / 2 = (1 - \nu)^2 / 2 \\ c_k(\nu) &= \langle \cos(\pi k \cdot), \Phi_\nu \rangle = \frac{2}{\pi^2 k^2} \left( \cos(\pi k) - \cos(\pi k \nu) \right) \\ b_k(\nu) &= \langle \sin(\pi k \cdot), \Phi_\nu \rangle = 0. \end{cases} \tag{2.2}$$

We insist here on the fact that the procedure applies for any other integrated functional $\int \Phi(f)$ with known continuous $\Phi$ when periodized. The values of the Fourier coefficients $c_k$ will change, but they will still be bounded by a quantity of order $k^{-2}$ for large $k$ and all the proofs work out the same way.

Let us discuss on the class constraints in (2.1). On the one hand, we assume densities $f$ to be uniformly bounded by some constant $\rho < 1$. More generally, we could have considered a class $\mathcal{F}(\mathbb{K}, R, m^*)$ (for $R > 0$ fixed) of probability density functions $f : \mathbb{K} \to \mathbb{R}$ such that $f(t) \geq m^* > 0$ for all $t \in \mathbb{K}$ and such that there exists some $0 < \rho < R$ and $\|f\|_\infty \leq \rho < R$. Then, the excess mass is defined via the functional $\Phi_\nu : [-R, R] \to [0, R]$ for any $0 \leq \nu \leq R$. For this functional we consider the rescaled Fourier basis on $[-R, R]$ and the corresponding coefficients are

$$c_0(\nu) = \frac{(R - \nu)^2}{2}, \quad c_k(\nu) = \frac{2R^2}{\pi^2 k^2}\left(\cos(\pi k) - \cos(\pi k \frac{\nu}{R})\right).$$

Therefore, without loss of generality we consider $R = 1$. On the other hand, we ask that the underlying density to be bounded from below away from 0. This is a classical assumption in the density model. Indeed, the variance of density estimators are proportional to $f$ and it cannot be controlled without such an assumption.

## 2.2. Estimation of the density $f$

We need now a nonparametric estimator of the density $f$. We can use any method and tune the smoothing parameter similarly. We chose the wavelet estimator of $f$ in order to deal easier with higher dimensions and to general functions in Besov classes. For this purpose, let us be given a pair of scaling function $\phi$ and associated wavelet function $\psi$. We assume that these functions are compactly supported (of support $[0, 2M]$); they can be of class $\mathcal{C}^r$ with $r$ as large as desired, see for example the Daubechies's wavelets, (6). With tensorial product, one can construct a multivariate scaling function and $2^d - 1$ associated wavelets always denoted by $\{\phi, \psi^\epsilon\}_{\epsilon \in \{1, \ldots, 2^d - 1\}}$, see (18). In the sequel, for any function $g \in L^2(\mathbb{R}^d), l \in \mathbb{Z}^{2d}, j \in \mathbb{N}$, we use the notation $g_{j,l}(.) = 2^{jd/2} g(2^j . - l)$. For a given $j \in \mathbb{N}$, the set $\{\phi_{j,l_1}, \psi^\epsilon_{j',l_2}, j' \geq j, (l_1, l_2) \in \mathbb{Z}^{2d}, \epsilon \in \{1, \ldots 2^d - 1\}\}$ is an orthonormal basis of $L^2(\mathbb{R}^d)$ and one can write, with the usual notations for the projections

$$\forall g \in L^2(\mathbb{R}^d), \quad \forall j \geq 0, \quad g = E_j g + D_j g \tag{2.3}$$

where

$$E_j g = \sum_{l \in Z^d} \alpha_{j,l} \phi_{j,l} \ \text{and} \ D_j g = \sum_{j' \geq j} \sum_{l \in Z^d} \sum_{1 \leq \epsilon \leq 2^d - 1} \beta^\epsilon_{j',l} \psi^\epsilon_{j',l}.$$

We omit the spaces where the indices are varying: $j$, $j\prime$ are always integers and $l$ is always a $d$−dimensional index. Denote $\lfloor . \rfloor$ the integer value and define $j_0 = \lfloor j_0' \rfloor$ and $j_\infty = \lfloor j_\infty' \rfloor$ where $j_0'$ and $j_\infty'$ are such that

$$2^{j_0'} = \log n \quad \text{and} \quad 2^{j_\infty'} = \left(\frac{n}{\log n}\right)^{1/d}. \tag{2.4}$$

Taking advantage of the decomposition (2.3), we propose to estimate $f$ by its wavelet estimate at the level $j$ varying between $j_0$ and $j_\infty$

$$\hat{f}_j(t) \stackrel{\text{def}}{=} \sum_l \hat{\alpha}_{j,l} \phi_{j,l}(t) \tag{2.5}$$

where the empirical coefficients are defined for any integer $j$ varying between $j_0$ and $j_\infty$ and for any integer $l \in \{2^j A_1 - 2M, 2^j B_1\} \times \ldots \times \{2^j A_d - 2M, 2^j B_d\}$ (we denoted $\mathbb{K} = [A_1, B_1] \times \ldots \times [A_d, B_d]$),

$$\hat{\alpha}_{j,l} = \frac{1}{n} \sum_{i=1}^n \phi_{j,l}(X_i).$$

Put

$$\lambda_j^2(t) = V(\hat{f}_j(t)) = \frac{1}{n} \sum_{l_1, l_2} \left(\int \phi_{j,l_1} \phi_{j,l_2} f - \int \phi_{j,l_1} f \int \phi_{j,l_2} f\right) \phi_{j,l_1}(t) \phi_{j,l_2}(t)$$

and observe that

$$\lambda_j^2(t) \leq \left((2M)^{2d} \|\phi\|_\infty^2 \|\phi\|_2^2 \|f\|_\infty\right) \frac{2^{jd}}{n}.$$

Using (2.1), we bound the constant in the right term by $\gamma = (2M)^{2d} \|\phi\|_\infty^2$. Moreover, we need to bound from below the variance $\lambda_j(t)$. Therefore, we choose a wavelet such that there exists $m > 0$ satisfying the assumption

$$\forall j = j_0 \ldots, j_\infty, \ \forall t \in \mathbb{K}, \ \exists l, \ |\phi(2^j t - l)| > m \tag{2.6}$$

where $j_0, j_\infty$ are defined in (2.4).

### 2.3. Plug-in

A candidate to estimate $A_N \Phi_\nu(f(t))$ could be $c_0(\nu) + \sum_{k=1}^N c_k(\nu) \cos(\pi k \hat{f}_j(t))$. Following (16), this estimator has too large a bias and we decrease this bias by considering the following modification. We estimate $A_N \Phi_\nu(f(t))$ by

$$A_{N,j}(t) = c_0(\nu) + \sum_{k=1}^N c_k(\nu) \exp(\pi^2 k^2 \lambda_j(t)^2/2) \cos(\pi k \hat{f}_j(t)).$$

Since the variance of the estimate of the density $\lambda_j(t)$ is unknown, we replace it with an estimate based on the empirical moments

$$\widehat{\lambda_j^2}(t) = \tag{2.7}$$

$$\frac{1}{n} \sum_{l_1, l_2} \left[\frac{1}{n} \sum_{i=1}^n \phi_{j,l_1}(X_i) \phi_{j,l_2}(X_i) - \frac{1}{n} \sum_{i=1}^n \phi_{j,l_1}(X_i) \frac{1}{n} \sum_{i=1}^n \phi_{j,l_2}(X_i)\right] \phi_{j,l_1}(t) \phi_{j,l_2}(t).$$

Notice that there exists some constant $c > 0$ such that $\widehat{\lambda_j^2}(t) \leq c\,2^{2jd}n^{-1}$ which could be much larger than $\lambda_j^2(t)$. We decide then to truncate $\widehat{\lambda_j^2}(t)$ at $\gamma 2^{jd}n^{-1}$ (for $\gamma = (2M)^{2d}\|\phi\|_\infty^2$) which is the upper bound for $\|\lambda_j^2(\cdot)\|_\infty$.

### 2.4. Estimator of the excess mass

Finally, we propose to estimate $\mathcal{E}(\nu)$ by

$$
\begin{aligned}
\hat{\mathcal{E}}(\nu) &= \int_{\mathbb{K}} \widehat{A_{N,j}}(t)dt \qquad\qquad\qquad\qquad\qquad (2.8)\\
&= \sum_{k=0}^{N} c_k(\nu) \int_{\mathbb{K}} \exp\left(\frac{\pi^2 k^2}{2} \min\left\{\widehat{\lambda_j^2}(t), \gamma\frac{2^{jd}}{n}\right\}\right) \cos\left(\pi k \hat{f}_j(t)\right) dt
\end{aligned}
$$

for $\widehat{\lambda_j^2}(t)$ defined in (2.7) and

$$
c_0(\nu) = (1-\nu)^2/2,\; c_k(\nu) = 2(\pi k)^{-2}(\cos(\pi k) - \cos(\pi k\nu)),\; \gamma = (2M)^{2d}\|\phi\|_\infty^2\,.
$$

## 3. Upper bounds and convergence properties

In Proposition 3.1, we give bounds from above for the expected errors of the estimation procedure of the functional $\mathcal{E}(\nu)$. This bound is depending on the parameters of estimation $j$ and $N$ and on the wavelet approximation error of $f$. Next, we determine the optimal parameters $j$ and $N$ to balance the terms appearing in the upper bound of Proposition 3.1, under the additional smoothness assumption on the unknown function $f$. In Theorem 3.1, an upper bound of order $(n \log n)^{-s/(2s+d)}$ for our estimation problem is found. In the gaussian white model nearly minimax lower bounds of order $(n \log n)^{-s/(2s+d)}(\log n)^{-1/(2s+1)}$ ($d = 1$, at a log factor) were established by (3). They improved the techniques used by (16) who found lower bounds of order $(n \log n)^{-s/(2s+d)}(\log n)^{-1}$ ($d = 1$) for estimating the $\mathbb{L}_1$ norm in the Gaussian white noise model, but there is still a gap between upper and lower bounds.

### 3.1. Expansion of the estimation error

The following bound for the mean absolute error of estimation of the excess mass holds.

**Proposition 3.1.** *Let $j$ be an integer between $j_0$ and $j_\infty$ and $N$ a positive integer. Assume that $f \in \mathcal{F}(\mathbb{K}, m^*)$. Choosing $\phi$ such that (2.6) holds for some*

$m > 0$, we get

$$
\begin{aligned}
E_f\left(d(\hat{\mathcal{E}}, \mathcal{E})\right) \;\leq\; & C_4 \frac{1}{N} + \|D_j f\|_1 \\
& + \left[ C_1\, N \left( \frac{2^{jd}}{n} \right)^{3/2} (\log n)^{1/2} + C_2\, \frac{1}{n^{1/2}} \log N \right. \\
& \left. + C_3 \frac{1}{N} \left( \frac{2^{jd}}{n} \right)^{1/2} \right] \exp\left( \frac{\pi^2 \gamma}{2} N^2 \frac{2^{jd}}{n} \right)
\end{aligned}
$$

*for $d$ denoting either i) the point wise difference, i.e. $d(g,h) = |g(\nu) - h(\nu)|$ for a given $\nu > 0$, ii) the sup-norm, or iii) the normalized $\mathbb{L}_2$−norm, i.e. $d(g,h) = \|g - h\|_2 / |\mathbb{K}|$ with $|\mathbb{K}|$ denoting the Lebesgue measure of the set $\mathbb{K}$ and*

$$
C_1 = 2|\mathbb{K}|, \quad C_2 = (4M)^{d/2} 4\pi^{-2} (2M)^d \|\phi\|_\infty, \quad C_3 = C_4 = 4\pi^{-2} |\mathbb{K}|.
$$

### 3.2. Upper bound for the estimation error

Let us now tune the parameters $N$ and $j$ in an optimal way. We denote, for fixed $m^* > 0$,

$$
\mathcal{F}(m^*) = \bigcup_{\mathbb{K}} \{ \mathcal{F}(\mathbb{K}, m^*) : |\mathbb{K}| \leq D \},
$$

for some fixed constant $D > 0$. We assume a Besov type smoothness condition for $f$ related to the wavelet expansion of the density $f$. More precisely, let $p, q \geq 1$, $s > 0$ and $L > 0$. The Besov bodies are characterized in term of wavelet coefficients as follows

$$
f \in b_{p,q}^s(L) \Leftrightarrow \|\alpha_{0,\cdot}\|_p + \left( \sum_{j \geq 0} \left[ 2^{j(s + \frac{d}{2} - \frac{d}{p})} \|\beta_{j \cdot}\|_p \right]^q \right)^{1/q} \leq L . \tag{3.1}
$$

Note that, for a given $r−$smooth wavelet with $r > s$, the Besov norm of a function $f$ is equivalent to the sequence norm of the wavelet coefficients of the function and then the concepts of Besov body and of Besov spaces are equivalent. For further details on the Besov spaces and their links with the wavelet analysis, see for instance (10). We derive from the smoothness assumption (3.1) on the unknown function $f$ the following bound for the bias term

$$
\exists \{\epsilon_j\}_{j \in \mathbb{N}} \in l_q, \quad \|D_j f\|_1 \leq \epsilon_j\, 2^{-js}.
$$

We choose the integer $j$ depending on $N$ such that $\lfloor 2^{js} \rfloor = N$ in order to balance the bias and the approximation error. We replace $j$ and minimize next the variance terms

$$
\frac{C_5}{N} + \left( C_1 N \left( \frac{N^{d/s}}{n} \right)^{3/2} \sqrt{\log n} + C_2 \frac{\log(N)}{\sqrt{n}} \right.
$$
$$
\left. + C_3 N^{-1} \left( \frac{N^{d/s}}{n} \right)^{1/2} \right) \exp\left( \frac{\pi^2 \gamma}{2} \frac{N^{(2s+d)/s}}{n} \right).
$$

We take $N = \lfloor (C_0 n \log n)^{s/(2s+d)} \rfloor$, with a constant $C_0 > 0$ such that $C_0 \pi^2 \gamma/2 < \min\{s, d/2\}/(2s+d)$. In this way, the exponential term in the variance term becomes a polynomial term smaller than the bias and the approximation term. The latter terms are of the same order: $(n \log n)^{-s/(2s+d)}$. Note that the variance term does not drive the rate. The following theorem is then proved.

**Theorem 3.1.** *Let $s > 0$, $1 \le p \le \infty$, $1 \le q \le \infty$, $L$, $D$, $m^* > 0$ and $0 < \rho < 1$. Let us suppose that $f$ belongs to $\mathcal{F}(m^*) \cap b^s_{p,q}(L)$ and assume there exists $m$ positive such that the technical assumption (2.6) holds. Let $\hat{\mathcal{E}}^*(\cdot)$ be the estimate of $\mathcal{E}(\cdot)$ defined by (2.5)-(2.8) for the following choice of estimation parameters*

$$j^* = \lfloor j^{*\prime} \rfloor, \ 2^{j^{*\prime}} = (n \log n)^{\frac{1}{2s+d}}, \ N^* = \lfloor (C_0 n \log n)^{s/(2s+d)} \rfloor$$

*where $C_0 > 0$ is a constant smaller than $\min\{2s, d\} \cdot \left( \pi^2 (2M)^{2d} \|\phi\|_\infty^2 (2s+d) \right)^{-1}$. Then*

$$\varlimsup_{n \to \infty} \sup_{f \in \mathcal{F}(m^*) \cap b^s_{p,q}(L)} (n \log n)^{\frac{s}{2s+d}} E_f \left( d(\hat{\mathcal{E}}^*, \mathcal{E}) \right) \ \le \ C,$$

*for $d$ denoting either i) the point wise difference, i.e. $d(g, h) = |g(\nu) - h(\nu)|$ for a given $\nu > 0$, ii) the sup-norm, or iii) the normalized $\mathbb{L}_2-$norm, i.e. $d(g, h) = \|g - h\|_2 / |\mathbb{K}|$ with $|\mathbb{K}|$ denoting the Lebesgue measure of the set $\mathbb{K}$ and the constant $C > 0$ depends on $s$, $L$, $D$, $d$ and $\phi$.*

As we already mentioned, the theorem is still valid if we estimate any other integrated functional of the type $\int \Phi(f)$ with $\Phi$ continuous not differentiable.

## 4. Numerical results

First, we describe the implementation of our estimation procedure $\hat{\mathcal{E}}^*(\nu)$ at level $\nu > 0$. In order to compare, we also implement a plug-in procedure $\hat{\mathcal{E}}^{PI}(\nu)$ defined as follows

$$\hat{\mathcal{E}}^{PI}(\nu) = \int_{\hat{f}_n(x) - \nu > 0} (\hat{f}_n(x) - \nu) dx.$$

where $\hat{f}_n$ is a density estimator. Let us recall that the best rate achievable by this procedure on the Besov balls (for the same loss functions as in Theorem 3.1) is the usual nonparametric rate $n^{-\frac{s}{2s+d}}$ obtained for the tuning parameter of order $n^{\frac{1}{2s+d}}$.

We compare several error measurements: in the sequel, $E_2^*$ and $E_\infty^*$ (respectively $E_2^{PI}$ and $E_\infty^{PI}$) denote the integrated squared error and the sup-norm due to our estimator $\hat{\mathcal{E}}^*$ (respectively due to the plug-in estimator $\hat{\mathcal{E}}^{PI}$). Moreover, the probability $p_2 = P_f(E_2^* < E_2^{PI})$ that the error of our procedure $\hat{\mathcal{E}}^*$ be smaller than the corresponding error of $\hat{\mathcal{E}}^{PI}$ is a good indicator of the performances of our procedure with respect to the plug-in procedure. Similarly, we consider $p_\infty = P_f(E_\infty^* < E_\infty^{PI})$

In the first part, we explain the automatic algorithm: observe that it is slightly different than the procedure described in the theoretical section (with respect to adaptation for instance). Next, we give a short summary of our simulation results: we try a lot of densities and we present here the most representative and relevant examples.

### *4.1. Algorithm*

The simulations are performed with the free software R V2.4. For the plug-in procedure, the estimator $\hat{f}_n$ is computed with the data driven procedure called *density()* provided by R. This kernel procedure of density estimation determines automatically the smoothing parameter $h^{**}$ that we use for the plug-in procedure. Since the theoretical optimal index for the plug-in procedure is $n^{-1/(d+2s)}$, we deduce $\hat{s}$ from $h^{**}$. Then we modify the smoothing index introducing the logarithmic term as indicated in Theorem 3.1. The parameters used when our procedure is computed are given by

$$\hat{N} = \lfloor (C_0 n \log n)^{\frac{\hat{s}}{d+2\hat{s}}} \rfloor, \ \hat{h} = (n \log n)^{-\frac{1}{d+2\hat{s}}}, \ C_0 = d.$$

We emphasize that the procedure *density()* is again used for our ouwn procedure but with the smoothing index modified as prescribed in Theorem 3.1. A bootstrap procedure of 100 replications is introduced to estimate the expected value $E_f(\hat{f}_n(x))$ and the variance $\hat{\lambda}^2 = V_f(\hat{f}_n(x))$ of $\hat{f}_n(x)$. As the bootstrap procedure gives a very accurate estimator of the variance, the truncation in the exponential term is actually useless for practical purposes and stands in formula (2.8) only for technical reasons in the proof. It is then sufficient to describe the estimator in (2.8) as

$$\sum_{k=0}^{\hat{N}} c_k(\nu) \int_{\mathbb{K}} \exp\left(\frac{\pi^2 k^2}{2} \hat{\lambda}^2(x)\right) \cos\left(\pi k \hat{f}_n(x)\right) dx.$$

As explained in the introduction, the exponential factor is a correction of the bias introduced when $\hat{f}$ is plugged-into the cosine function. Therefore, we suggest to compute the estimator as

$$\hat{\mathcal{E}}^*(\nu) = \sum_{k=0}^{N} c_k(\nu) \int_{\mathbb{K}} \cos\left(\pi k E_f(\hat{f}_n(x))\right) dx,$$

where $E_f(\hat{f}_n(x))$ is very well recovered by a bootstrap estimation procedure. In practice, both methods give the same results, but the second formula is computed significantly faster than the first.

A sequence $0 = \nu_1, \ldots, \nu_{100} = 1$ is considered. The empirical errors denoted $\widetilde{E}_2^*$, $\widetilde{E}_2^{PI}$ and $\widetilde{E}_\infty^*$, $\widetilde{E}_\infty^{PI}$ are computed via $K = 20$ Monte Carlo simulations. We denote $\widetilde{p}_2$ and $\widetilde{p}_\infty$ the frequencies of success of our procedure.

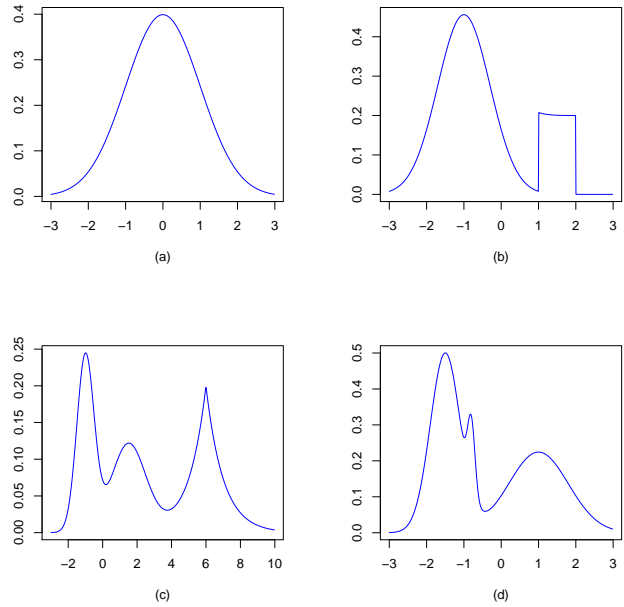

Fig 2. *Set of studied densities. (a): standard gaussian; (b): mixture of gaussian and uniform; (c): mixture of 2 gaussian and laplace; (d):mixture of gaussian with isolated spoke.*

## 4.2. Univariate densities

We consider some example of probability densities which are mixtures of gaussian, uniform and Laplace laws, see Figure 2. The density (a)

$$f(x) = f_{\mathcal{N}(0,1)}(x)$$

is the gaussian density: this is the most standard example and it is very popular in practical studies. The density (b)

$$f(x) = 0.8 \cdot f_{\mathcal{N}(-1,0.7)}(x) + 0.2 \cdot f_{\mathcal{U}(1,2)}(x)$$

is a mixture of a gaussian density and a uniform density. Remark that the uniform density is not continuous and this allows us to study the robustness of our procedure. Moreover, the gaussian part is very smooth: the mixture density is then difficult to estimate because there is a conflict about a global choice of the bandwidth. The density (c)

$$f(x) = 0.3 \cdot f_{\mathcal{N}(-1,0.5)}(x) + 0.3 \cdot f_{\mathcal{N}(1.5,1)}(x) + 0.4 \cdot f_{\mathcal{L}(6)}(x)$$

is a mixture of a gaussian density and Laplace density. Since the Laplace density is not differentiable at its mode, we study again the same phenomenon: the *Univariate densities. Comparison of $\hat{\mathcal{E}}^*$ and $\hat{\mathcal{E}}^{PI}$ in mean integrated squared error and in mean error of the sup-norm, over $K = 20$ Monte-Carlo simulations, for various sizes of samples $n = 100, 1000, 10000$.*

| $f$ | $n$ | $\widetilde{E}_2^{PI}$ | $\widetilde{E}_2^*$ | $\widetilde{E}_2^{PI}/\widetilde{E}_2^*$ | $\widetilde{p}_2$ | $\widetilde{E}_\infty^{PI}$ | $\widetilde{E}_\infty^*$ | $\widetilde{E}_\infty^{PI}/\widetilde{E}_\infty^*$ | $\widetilde{p}_\infty$ |
|---|---|---|---|---|---|---|---|---|---|
| a | 100 | 0.00504 | 0.00542 | 0.93 | 0.45 | 0.04450 | 0.04855 | 0.92 | 0.45 |
| a | 1000 | 0.00079 | 0.00066 | 1.19 | 0.70 | 0.01937 | 0.01765 | 1.10 | 0.75 |
| a | 10000 | 0.00008 | 0.00006 | 1.32 | 0.55 | 0.00602 | 0.00590 | 1.02 | 0.60 |
| b | 100 | 0.00354 | 0.00533 | 0.66 | 0.20 | 0.03742 | 0.04989 | 0.75 | 0.30 |
| b | 1000 | 0.00147 | 0.00086 | 1.71 | 0.90 | 0.03217 | 0.02133 | 1.51 | 0.95 |
| b | 10000 | 0.00170 | 0.00083 | 2.06 | 1.00 | 0.03645 | 0.02445 | 1.49 | 1.00 |
| c | 100 | 0.00520 | 0.01027 | 0.51 | 0.15 | 0.04132 | 0.04924 | 0.84 | 0.30 |
| c | 1000 | 0.00077 | 0.00036 | 2.17 | 0.80 | 0.01745 | 0.01274 | 1.37 | 0.80 |
| c | 10000 | 0.00075 | 0.00021 | 3.64 | 1.00 | 0.01714 | 0.00938 | 1.83 | 1.00 |
| d | 100 | 0.03271 | 0.01857 | 1.76 | 1.00 | 0.11473 | 0.08293 | 1.38 | 1.00 |
| d | 1000 | 0.00975 | 0.00346 | 2.82 | 1.00 | 0.05985 | 0.03606 | 1.66 | 1.00 |
| d | 10000 | 0.00248 | 0.00063 | 3.91 | 1.00 | 0.02975 | 0.01525 | 1.95 | 1.00 |

smoothing indices can not be at the same time globally designed and everywhere optimal. The last density (d)

$$f(x) = 0.5 \cdot f_{\mathcal{N}(-1.5, 0.4)}(x) + 0.05 \cdot f_{\mathcal{N}(-0.8, 0.1)}(x) + 0.45 \cdot f_{\mathcal{N}(1, 0.8)}(x)$$

is a mixture of three gaussian densities with isolated peaks and different variances. Density (d) is a case where the smoothing indices of the estimation procedures have to be space-dependant in view to capture the small sharp peak.

One challenge is to check whether our procedure overcomes all the enumerated difficulties for the estimation of the density $f$. The results are presented in Table 1.

First, we note that our procedure is becoming more accurate when the size of the sample increases. It seems that our method is relatively complicated and need enough data to be powerful. In the opposite, the naive plug-in method is a robust procedure which is not so bad when few data are available: when $n = 100$, the frequencies of success of the plug-in method with respect to our procedure is $1 - \tilde{p} = 0.55, 0.80, 0.85$ for the density (a), the density (b) and the density (c). But when $n$ is larger, our method is more successful: $\tilde{p} = 0.90, 0.80$ for the density (b), the density (c).

When the densities become more and more complex (by complex, we mean an increase of the number of modes or irregularities in the density), our procedure is much more relevant than the plug-in procedure. For large samples, $n = 10000$, we observe a benefit of 106% for a mixture of gaussian and uniform densities, a benefit of 264% for a mixture of gaussian and Laplace and a benefit of 291% for a mixture of densities with a small isolated peak. In parallel, we observe that, for all the Monte Carlo simulations, our estimator $\hat{\mathcal{E}}^*$ is systematically better than $\hat{\mathcal{E}}^{PI}$, with a rate $\tilde{p}_2 = \tilde{p}_\infty = 1$.

Observe that for the density (d), our method is better than the plug-in estimator for any sample size. It seems that the change of the smoothing parameter adding an extra logarithmic term is crucial to kill a great part of the bias term.

### 4.3. Bivariate densities

In this part, we focus on gaussian and uniform densities. Let us denote

$$\mathcal{N}((EX, EY), (\sqrt{V(X)}, \sqrt{V(Y)}, \rho_{XY}))$$

the bivariate gaussian density of $(X, Y)$. The studied densities are plotted in Figure 3. The density (A) is the standard one

$$f = f_{\mathcal{N}((0,0),(1,1,0))}.$$

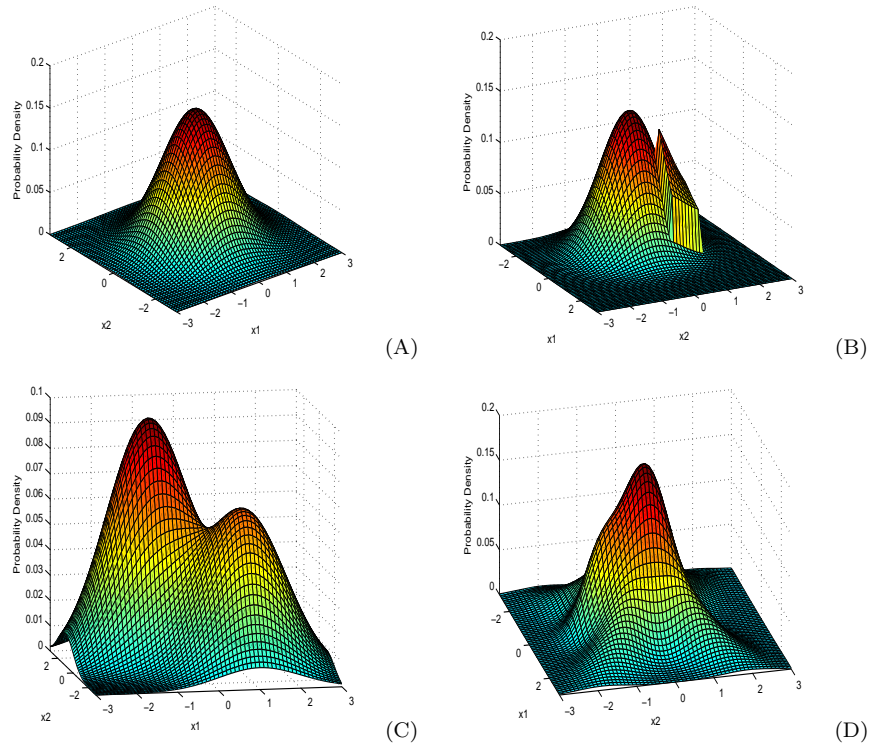

FIG 3. *Set of studied 2D densities. (A): 2D Gaussian. (B): mixture of 2D gaussian and uniform. (C): mixture of two 2D gaussian. (D): Mixture of three 2D gaussian.*

Bivariate densities. Comparison of $\tilde{\mathcal{E}}^*$ and $\tilde{\mathcal{E}}^{PI}$ in mean integrated squared error and in mean error of the sup-norm, over $K = 20$ Monte-Carlo simulations, for various sizes of samples $n = 400, 1000, 10000$.

| $f$ | $n$ | $\widetilde{E}_2^{PI}$ | $\widetilde{E}_2^*$ | $\widetilde{E}_2^{PI}/\widetilde{E}_2^*$ | $\widetilde{p}_2$ | $\widetilde{E}_\infty^{PI}$ | $\widetilde{E}_\infty^*$ | $\widetilde{E}_\infty/\widetilde{E}_\infty^*$ | $\widetilde{p}_\infty$ |
|---|---|---|---|---|---|---|---|---|---|
| A | 400 | 0.01685 | 0.00870 | 1.94 | 0.95 | 0.07435 | 0.05675 | 1.31 | 0.95 |
| A | 1000 | 0.00948 | 0.00394 | 2.41 | 1.00 | 0.05445 | 0.03525 | 1.54 | 1.00 |
| A | 10000 | 0.00635 | 0.00263 | 2.41 | 1.00 | 0.04524 | 0.02867 | 1.58 | 1.00 |
| B | 400 | 0.10397 | 0.04628 | 2.25 | 1.00 | 0.17667 | 0.12818 | 1.38 | 1.00 |
| B | 1000 | 0.07460 | 0.02943 | 2.53 | 1.00 | 0.15398 | 0.10660 | 1.44 | 1.00 |
| B | 10000 | 0.04747 | 0.01985 | 2.39 | 1.00 | 0.12894 | 0.08901 | 1.45 | 1.00 |
| C | 400 | 0.01184 | 0.00555 | 2.13 | 1.00 | 0.06442 | 0.04609 | 1.40 | 1.00 |
| C | 1000 | 0.00906 | 0.00432 | 2.09 | 1.00 | 0.05462 | 0.03717 | 1.47 | 1.00 |
| C | 10000 | 0.00379 | 0.00134 | 2.83 | 1.00 | 0.03566 | 0.02081 | 1.71 | 1.00 |
| D | 400 | 0.26965 | 0.26187 | 1.03 | 0.55 | 0.34953 | 0.34350 | 1.02 | 0.55 |
| D | 1000 | 0.25655 | 0.24609 | 1.04 | 0.85 | 0.33525 | 0.32628 | 1.03 | 0.85 |
| D | 10000 | 0.23705 | 0.22277 | 1.06 | 1.00 | 0.31686 | 0.30497 | 1.04 | 1.00 |

The density (B) is a mixture of the Gaussian density and the uniform density

$$f = 0.6 \cdot f_{\mathcal{N}((-1,0),(0.7,0.7,0))} + 0.4 \cdot f_{\mathcal{U}([0.5,1.5] \times [-0.5,0.5])}$$

and Density (C)

$$f = 0.8 \cdot f_{\mathcal{N}((-0.5,0.5),(1,1,0))} + 0.2 \cdot f_{\mathcal{N}((0.4,-0.4),(1,1,0))}$$

is a mixture of two gaussian densities. Last, the density (D)

$$\begin{aligned} f &= 0.45 \cdot f_{\mathcal{N}((0,0),(1.5,1,0.95))} + 0.45 \cdot f_{\mathcal{N}((0,0),(1.5,1,-0.95))} \\ &\quad + 0.10 \cdot f_{\mathcal{N}((0,-1.2),(0.2,0.2,0))} \end{aligned}$$

is a mixture of three gaussian densities presenting an isolated spot.

Table 2 presents results for different sample sizes $n = 400, 1000, 10000$. As in the one dimensional case, we observe that the improvement increases with $n$ and with the complexity of the underlying theoretical density increases. We observe that our procedure is always better than the plug-in procedure. The frequencies of success of our procedure are very high: $\tilde{p}_2 = \tilde{p}_\infty = 0.55$ for the densities (A), (B), (C). Even in the case of the density (D) where the results are mitigated, the worse result is $\tilde{p}_2 = \tilde{p}_\infty = 0.55$ when $n = 400$ (which is very small for $2$−dimensional non parametric estimation problems).

The improvements are more remarkable in bivariate case (but for the density (D)). We may say that the gain of a logarithmic factor in the rate of our estimator gives a significant compensation for the curse of dimensionality. For the standard gaussian density, we observe a high improvement of our estimation procedure compared to the plug-in estimator as the dimension increases: for $n = 10000$, the empirical mean squared error has improved from 32% in the case of the density (a) to 141% in the case of the density (A). When empirical sup-norm is considered, the improvements are not so extraordinary but there are significant: from 2% to 58%.

We think that the mediocrity of the results in the case of the density (D) could be corrected by a more accurate determination of the smoothing indices.

## 5. Proofs

### 5.1. Proof of Proposition 3.1

In order to study the quadratic error, it is useful to note that the hypothesis (2.1) for some $\rho < 1$ implies that $\mathcal{E}(\nu)$ is zero if $\nu \geq 1$. Let us summarize again some notation

- $\Phi_\nu(t) = (|t| - \nu)\, 1_{|t| > \nu}$,
- $A_N \Phi_\nu(t) = c_0(\nu) + \sum_{k=1}^{N} c_k(\nu) \cos(\pi k t)$,
- $A_{N,j}(t, \nu) = c_0(\nu) + \sum_{k=1}^{N} c_k(\nu) \exp(\pi^2 k^2 \lambda_j(t)^2/2) \cos(\pi k \hat{f}_j(t))$,
- $\widehat{A_{N,j}}(t, \nu) = c_0(\nu) + \sum_{k=1}^{N} c_k(\nu) \exp\left(\frac{\pi^2 k^2}{2} \max\left\{\widehat{\lambda_j^2}(t), \gamma \frac{2^{jd}}{n}\right\}\right) \cos(\pi k \hat{f}_j(t))$,
- the unknown density $f$ writes on the wavelet basis $f = E_j f + D_j f$ where

$$E_j f = \sum_l \alpha_{j,l} \phi_{j,l} \quad \text{and} \quad D_j f = \sum_{j' \geq j} \sum_l \sum_\varepsilon \beta_{j',l}^\varepsilon \psi_{j',l}^\varepsilon.$$

Then $\hat{\mathcal{E}}(\nu) = \int_{\mathbb{K}} \widehat{A_{N,j}}(t, \nu) dt$. We have the following expansion

$$\hat{\mathcal{E}}(\nu) - \mathcal{E}(\nu) = S_1(\nu) + S_2(\nu) + B_2(\nu) + A(\nu) + B_1(\nu)$$
$$= \int_{\mathbb{K}} \left[\widehat{A_{N,j}}(t, \nu) - A_{N,j}(t, \nu)\right] dt + \int_{\mathbb{K}} \left[A_{N,j}(t, \nu) - E(A_{N,j}(t, \nu))\right] dt$$
$$+ \int_{\mathbb{K}} \left[E(A_{N,j}(t, \nu)) - A_N \Phi_\nu(E_j f(t))\right] dt$$
$$+ \int_{\mathbb{K}} \left[A_N \Phi_\nu(E_j f(t)) - \Phi_\nu(E_j f(t))\right] dt$$
$$+ \int_{\mathbb{K}} \left[\Phi_\nu(E_j f(t)) - \Phi_\nu(f(t))\right] dt$$

where $S_1(\nu), S_2(\nu)$ are stochastic terms, $B_2(\nu)$ is a bias term due to the plug in, $A(\nu)$ is an approximation term and $B_1(\nu)$ is a bias term due to the estimation of the function of interest $f$. Proposition 3.1 is proved combining (5.1), (5.2), (5.3), (5.7) and (5.9).

*5.1.1. Bias term (due to the estimation).*

The bias term $B_1(\nu)$ is bounded using the fact that $|a(t)_+ - b(t)_+| \leq |a(t) - b(t)|$

$$|B_1(\nu)| \leq \int_{\mathbb{K}} |E_j f(t) - f(t)| \, dt \leq \|D_j f\|_1. \tag{5.1}$$

Note that the same bound holds for $\|B_1\|_\infty$ and for $\|B_1\|_2/|\mathbb{K}|$.

*5.1.2. Approximation term.*

Using the values of the Fourier coefficients given in (2.2), we have the following approximation for any $N$,

$$\forall u \in [-1, 1], \ |\Phi_\nu(u) - A_N \Phi_\nu(u)| \leq \frac{4}{\pi^2} \sum_{k > N} \frac{1}{k^2} \leq \frac{4}{\pi^2 N}$$

implying that

$$|A(\nu)| \ \leq \ \int_{\mathbb{K}} |A_N \Phi_\nu(E_j f(t)) - \Phi_\nu(E_j f(t))| \, dt \leq \frac{4|\mathbb{K}|}{\pi^2 N}. \tag{5.2}$$

Note that the same bound holds for $\|A\|_\infty$ and for $\|A\|_2/|\mathbb{K}|$.

*5.1.3. Bias term (due to the plug-in).*

First, we state the following lemma proved in the next section.

**Lemma 5.1.** *Let $N$ be a positive integer, $f$ belongs to $\mathcal{F}(\mathbb{K}, m^*)$ and $\hat{f}_j$ be the wavelet estimator constructed in (2.5). For $k = 1, \ldots, N$ and $j$ varying between $j_0$ and $j_\infty$, we have*

$$\forall t \in \mathbb{K}, \ \left| E\left[ e^{\pi^2 k^2 \lambda_j^2(t)/2} \cos(\pi k \hat{f}_j(t)) \right] - \cos(\pi k E_j f(t)) \right| \ \leq \ u_n \, e^{\pi^2 k^2 \lambda_j^2(t)/2}$$

*where*

$$u_n \ = \ \pi A \left( \frac{2^{jd}}{n} \right)^{1/2}.$$

*for an universal positive constant $A$.*

Applying Lemma 5.1, we get

$$
\begin{aligned}
|B_2(\nu)| \ &\leq \ \int_{\mathbb{K}} |E A_{N,j}(t, \nu) - A_N \Phi_\nu(E_j f(t))| \, dt \\
&\leq \ \sum_{k=1}^{N} |c_k(\nu)| \int_{\mathbb{K}} \left| e^{\pi^2 k^2 \lambda_j^2(t)/2} E(\cos(\pi k \hat{f}_j(t))) - \cos(\pi k E_j f(t)) \right| dt \\
&\leq \ |\mathbb{K}| \sum_{k=1}^{N} |c_k(\nu)| \sup_{t \in \mathbb{K}} \left( \sup_{k=1,\ldots,N} e^{\pi^2 k^2 \lambda_j^2(t)/2} u_n \right).
\end{aligned}
$$

Taking

$$\lambda_j^2(t) \leq \gamma \frac{2^{jd}}{n}, \ |c_k(\nu)| \leq \frac{4}{\pi^2 k^2}, \ u_n = \pi A \left( \frac{2^{jd/2}}{n} \right)^{1/2}$$

we obtain

$$|B_2(\nu)| \quad \leq \quad 4\pi^{-2} |\mathbb{K}| \frac{1}{N} \left( \frac{2^{jd}}{n} \right)^{1/2} \exp \left( \frac{\pi^2 \gamma}{2} N^2 \frac{2^{jd}}{n} \right). \tag{5.3}$$

Note that the same bound holds for $\|B_2\|_\infty$ and for $\|B_2\|_2 / |\mathbb{K}|$.

### 5.1.4. Stochastic term.

The wavelet estimator $\hat{f}_j(t)$ at point $t = (t_1, \ldots, t_d)$ is depending on the observations $X_i = (X_{1i}, \ldots, X_{di})$ such $|X_{pi} - t_i| \leq 2M \, 2^{-j}$ for any $p = 1, \ldots, d$. Therefore, $\hat{f}_j(t)$ and $\hat{f}_j(t')$ are independent as soon as there exists a direction $p$ such that $\|t_p - t'_p\| > 2M \, 2^{-j}$ and the same holds for any statistics $Z(t)$ and $Z(t')$ based on the observations. As in (16), adapted for $d$-dimensional setup:

$$E \left| \int_{\mathbb{K}} (Z(t) - EZ(t)) \, dt \right| \leq V^{1/2} \left( \int_{\mathbb{K}} Z(t) dt \right)$$

$$\leq \left( \int_{\mathbb{K}} \int_{\mathbb{K}} Cov \left( Z(t), Z(t') \right) dt dt' \right)^{1/2}$$

$$\leq \left( \int_{\mathbb{K}} \int_{\mathbb{K}} [V(Z(t)) \cdot V(Z(t'))]^{1/2} \prod_{p=1}^{d} I(|t_p - t'_p| \leq 2M2^{-j}) dt dt' \right)^{1/2}$$

$$\leq \left( \frac{1}{2} \int_{\mathbb{K}} \int_{\mathbb{K}} [V(Z(t)) + V(Z(t'))] \prod_{p=1}^{d} I(|t_p - t'_p| \leq 2M2^{-j}) dt dt' \right)^{1/2}$$

$$\leq (4M)^{d/2} 2^{-jd/2} \left( \int_{\mathbb{K}} V(Z(t)) \, dt \right)^{1/2}.$$

Denoting $Z_{j,k}(t) = \exp(\pi^2 k^2 \lambda_j(t)^2 / 2) \cos(\pi k \hat{f}_j(t))$, it follows

$$E(|S_2(\nu)|) \quad \leq \quad \sum_{k=1}^{N} c_k(\nu) E(| \int (Z_{j,k}(t) - E(Z_{j,k}(t))) dt |)$$

$$\leq \quad (4M)^{d/2} 2^{-jd/2} |\mathbb{K}| \sum_{k=1}^{N} |c_k(\nu)| \sup_{t \in \mathbb{K}} V^{1/2}(Z_{j,k}(t)) \tag{5.4}$$

$$E(\|S_2\|_\infty) \quad \leq \quad (4M)^{d/2} 2^{-jd/2} |\mathbb{K}| \sum_{k=1}^{N} \sup_{\nu} |c_k(\nu)| \sup_{t \in \mathbb{K}} V^{1/2}(Z_{j,k}(t)) \tag{5.5}$$

$$E(\|S_2\|_2) \quad \leq \quad (4M)^{d/2} 2^{-jd/2} |\mathbb{K}| \sum_{k=1}^{N} \|c_k(\cdot)\|_2 \sup_{t \in \mathbb{K}} V^{1/2}(Z_{j,k}(t)) \tag{5.6}$$

We state now the following lemma which is proved in the next section.

**Lemma 5.2.** *Let $N$ be a positive integer, $f$ belongs to $\mathcal{F}(\mathbb{K}, m^*)$ and $\hat{f}_j$ be the wavelet estimator constructed in (2.5). For $k = 1, \ldots, N$ and for the integer $j$ varying between $j_0$ and $j_\infty$*

$$\forall t \in \mathbb{K},\ V\left(e^{\pi^2 k^2 \lambda_j^2(t)/2} \cos(\pi k \hat{f}_j(t))\right) \leq \left(u_n + \pi^2 k^2 \lambda_j^2(t)\right) e^{\pi^2 k^2 \lambda_j^2(t)}$$

*where $u_n$ is given in Lemma 5.1.*

Direct application of Lemma 5.2 with the bound $\lambda_j(t)^2 \leq \gamma \left(\frac{2^{jd}}{n}\right)$ combined with (5.4) leads to

$$
\begin{aligned}
&E(|S_2(\nu)|) \\
\leq\ & (4M)^{d/2} 2^{-jd/2} |\mathbb{K}| \frac{4}{\pi^2} \left[\sum_{k=1}^N \frac{u_n^{1/2}}{k^2} e^{\pi^2 k^2 \lambda_j^2(t)/2} + \sum_{k=k_n}^N \frac{\pi \lambda_j(t)}{k} e^{\pi^2 k^2 \lambda_j^2(t)/2}\right] \\
\leq\ & (4M)^{d/2} 2^{-jd/2} |\mathbb{K}| \frac{4}{\pi^2} \left[u_n^{1/2} N^{-1} e^{\pi^2 N^2 \lambda_j^2(t)/2} + \pi \lambda_j(t) \log N\, e^{\pi^2 N^2 \lambda_j^2(t)/2}\right] \\
\leq\ & (4M)^{d/2} 2^{-jd/2} \frac{4\sqrt{\gamma}}{\pi} |\mathbb{K}| \left(\frac{2^{jd}}{n}\right)^{1/2} \log N \, \exp\left(\frac{\pi^2 \gamma}{2} N^2 \frac{2^{jd}}{n}\right).
\end{aligned}
$$

By (5.5) and (5.6), we obtain the same bound for $E(\|S_2\|_\infty)$ and $E(\|S_2\|_2)/|\mathbb{K}|$.

### 5.1.5. Stochastic term due to the estimation of the variance

Let $\widetilde{\lambda_j^2}(t) = \min\{\widehat{\lambda_j^2}(t), \gamma\, 2^{jd}/n\}$. We get

$$
\begin{aligned}
|S_1(\nu)| &\leq \int_{\mathbb{K}} \left|\widehat{A_{N,j}}(t, \nu) - A_{N,j}(t, \nu)\right| dt \\
&\leq \sum_{k=1}^N \left(|c_k(\nu)| \int_{\mathbb{K}} \left|\exp(\pi^2 k^2\, \widetilde{\lambda_j^2}(t)/2) - \exp(\pi^2 k^2\, \lambda_j^2(t)/2)\right| dt\right) \\
&\leq \sum_{k=1}^N \left(|c_k(\nu)| \frac{\pi^2 k^2}{2} \int_{\mathbb{K}} |\widetilde{\lambda_j^2}(t) - \lambda_j^2(t)| \exp(\pi^2 k^2\, \Lambda_{j,n}(t)/2) dt\right) \\
&\leq 2\sum_{k=1}^N \left(\int_{\mathbb{K}} |\widetilde{\lambda_j^2}(t) - \lambda_j^2(t)| \exp(\pi^2 k^2\, \Lambda_{j,n}(t)/2) dt\right) \\
&\leq 2N\, \exp\left(\frac{\gamma \pi^2}{2} N^2 \frac{2^{jd}}{n}\right) \int_{\mathbb{K}} |\widetilde{\lambda_j^2}(t) - \lambda_j^2(t)| dt \qquad (5.7)
\end{aligned}
$$

where $\Lambda_{j,n}(t)$ is an intermediate point between $\lambda_j^2(t)$ and $\widetilde{\lambda_j^2}(t)$ and therefore $|\Lambda_{j,n}(t)| \leq \gamma\, 2^{jd}/n$. We note that a rough bound like $|\widetilde{\lambda_j^2}(t) - \lambda_j^2(t)| \leq \gamma\, 2^{jd}/n$ is too large; thus, we take $\tau > 0$ and we split the expected value and for $q = 1, 2$, we have

$$
\begin{aligned}
E(|\widetilde{\lambda_j^2}(t) - \lambda_j^2(t)|) &\leq \tau + E(|\widetilde{\lambda_j^2}(t) - \lambda_j^2(t)| \, 1_{\{|\widetilde{\lambda_j^2}(t) - \lambda_j^2(t)| \geq \tau\}}) \\
&\leq \tau + \int_\tau^{\gamma \frac{2^{jd}}{n}} x \, d\tilde{P}(x)
\end{aligned}
$$

where $\tilde{P}$ is the probability associated with the variable $|\widetilde{\lambda_j^2}(t) - \lambda_j^2(t)|$. It follows that

$$
\begin{aligned}
E(|\widetilde{\lambda_j^2}(t) - \lambda_j^2(t)|) &\leq \tau + \gamma \frac{2^{jd}}{n} \int_\tau^{\gamma \frac{2^{jd}}{n}} d\tilde{P}(x) \\
&\leq \tau + \gamma \frac{2^{jd}}{n} P\left(|\widehat{\lambda_j^2}(t) - \lambda_j^2(t)| \geq \tau\right) \quad (5.8)
\end{aligned}
$$

since $|\widetilde{\lambda_j^2}(t) - \lambda_j^2(t)| \leq |\widehat{\lambda_j^2}(t) - \lambda_j^2(t)|$. We need an evaluation of the accuracy of the estimation of the variance of the estimate of $f$. The following lemma gives deviations for the error of estimation. For the proof, we refer to (25).

**Lemma 5.3.** *Let $f$ belong to $\mathcal{F}(m^*)$ and $t$ be in $\mathbb{K}$. Assume that there exists a positive constant $m$ such that (2.6) is satisfied. Then, for $\tau > 2^{jd} n^{-2}$, the estimator $\widehat{\lambda_j^2}(t)$ in (2.7) is such that*

$$
\begin{aligned}
P\left(|\widehat{\lambda_j^2}(t) - \lambda_j^2(t)| \geq \tau\right) &\leq c' \left[ \left(\frac{2^{jd}}{n}\right)^4 \frac{1}{\tau^2} \right. \\
&\left. + \exp\left\{-c\left(\left(\frac{2^{jd}}{n}\right)^{-3} \tau^2 \wedge \left(\frac{2^{jd}}{n}\right)^{-2} \tau\right)\right\} \right],
\end{aligned}
$$

*for some constants $c, c' > 0$.*

We use Lemma 5.3 with

$$
\tau = a \left(\frac{2^{jd}}{n}\right)^{3/2} (\log n)^{1/2}, \quad a > 0
$$

(which is larger than $\frac{2^{jd}}{n^2}$) to give an upper bound for (5.8)

$$
\begin{aligned}
&E(|\widetilde{\lambda_j^2}(t) - \lambda_j^2(t)|) \\
&\leq \left[ a \left(\frac{2^{jd}}{n}\right)^{3/2} (\log n)^{1/2} + \gamma c' a^{-2} \left(\frac{2^{jd}}{n}\right)^2 (\log n)^{-1} \right. \\
&\left. + \gamma c' \frac{2^{jd}}{n} \exp\left(-ca^2 \log n\right) + \gamma c' \frac{2^{jd}}{n} \exp\left(-ca\left(\frac{2^{jd}}{n \log n}\right)^{-1/2}\right) \right] \\
&\leq a \left(\frac{2^{jd}}{n}\right)^{3/2} (\log n)^{1/2}.
\end{aligned}
$$

since $j \leq j_\infty$ and as soon as $a \geq \max\{\sqrt{1/(2c)}, 1/(2c)\}$. From (5.7), we deduce

$$E(|S_1(\nu)|) \quad \leq \quad 2|\mathbb{K}|a\,N\,\exp\left(\frac{c\pi^2}{2}N^2\frac{2^{jd}}{n}\right)\left(\frac{2^{jd}}{n}\right)^{3/2}(\log n)^{1/2} \quad (5.9)$$

and the same bound is valid for $E(\|S_1\|_\infty)$. Observing that

$$E(\|S_1\|_2) \quad \leq \quad \sum_{k=1}^N \left(\|c_k(\cdot)\|_2 E \int_\mathbb{K} \left|\exp(\pi^2 k^2\,\widetilde{\lambda_j^2}(t)/2) - \exp(\pi^2 k^2\,\lambda_j^2(t)/2)\right| dt\right)$$

$$\leq \quad 2N\,\exp\left(\frac{\gamma\pi^2}{2}N^2\frac{2^{jd}}{n}\right)\,\int_\mathbb{K} E(|\widetilde{\lambda_j^2}(t) - \lambda_j^2(t)|)\,dt \quad (5.10)$$

and we obtain the same bound for $E(\|S_1\|_2)$.

### *5.2. Proofs of the lemmas*

#### *5.2.1. Proof of Lemma 5.1*

Let $t \in \mathbb{K}$ be fixed. Denote $K(k,j) = e^{\pi^2 k^2 \lambda_j^2/2}$. Recall that $\lambda_j^2(t) = V(\hat{f}_j(t))$, put

$$\chi_j = \frac{\hat{f}_j(t) - E_j f(t)}{\lambda_j}$$

and write

$$\cos(\pi k\hat{f}_j(t)) = \cos\left(\pi k E_j f(t) + \pi k\lambda_j\,\chi_j\right).$$

Expand using the formula of $\cos(a+b)$

$$\begin{aligned} z \quad &\overset{\text{def}}{=} \quad K(k,j)\,E(\cos(\pi k\hat{f}_j(t))) - \cos(\pi k E_j f(t)) \\ &= \quad [K(k,j)\,E\left(\cos\left(\pi k\lambda_j\,\chi_j\right)\right) - 1]\;\cos(\pi k E_j f(t)) \\ &\quad - K(k,j)\,E\left(\sin\left(\pi k\lambda_j\,\chi_j\right)\right)\sin(\pi k E_j f(t)). \end{aligned}$$

Observe that for $\chi$ a standard gaussian variable

$$K(k,j)\,E\left(\cos\left(\pi k\lambda_j\chi\right)\right) = 1 \quad \text{and} \quad K(k,j)\,E\left(\sin\left(\pi k\lambda_j\chi\right)\right) = 0.$$

We use an approximation for the law of $\chi_j$ as $n$ grows to $\infty$. Denote $F_{\mathcal{N}(0,1)}$ and $F_{\chi_j}$ the distribution functions of $\chi$ and $\chi_j$. It follows

$$\begin{aligned} z \quad = \quad &K(k,j)\cos(\pi k E_j f(t))\,\int \cos\left(\pi k\lambda_j\,x\right)(F_{\chi_j} - F_{\mathcal{N}(0,1)})(x)dx \\ &- K(k,j)\sin(\pi k E_j f(t))\,\int \sin\left(\pi k\lambda_j\,x\right)(F_{\chi_j} - F_{\mathcal{N}(0,1)})(x)dx \end{aligned}$$

Notice that

$$\chi_j = \frac{1}{\lambda_j(t)} \sum_l (\hat{\alpha}_{j,l} - \alpha_{j,l}) \phi_{j,l}(t) = \sum_{i=1}^n Z_{i,n}(t)$$

for

$$Z_{i,n}(t) = \frac{1}{n\,\lambda_j(t)} \left( \sum_l (\phi_{j,l}(X_i) - \alpha_{j,l}) \phi_{j,l}(t) \right).$$

Straightforward computations lead to

$$E(Z_{i,n}(t)^2) = \frac{1}{n}$$

and

$$
\begin{aligned}
E(|Z_{i,n}(t)|^3) &\leq \frac{2}{n^3 \lambda_j^3(t)} \sum_{l_1,l_2,l_3} \left| \left( \int \phi_{j,l_1} \phi_{j,l_2} \phi_{j,l_3} f \right) \phi_{j,l_1}(t) \phi_{j,l_2}(t) \phi_{j,l_3}(t) \right| \\
&\leq c\, \frac{2^{2jd}}{n^3 \lambda_j^3(t)}.
\end{aligned}
$$

for $c = (2M)^{3d} \|\phi\|_\infty^3 \|\phi\|_2 \|\phi\|_4^2 \|f\|_\infty$. Using the technical assumption (2.6) and the fact that there exists $m^* > 0$ such that $\inf_{t \in \mathbb{K}} f(t) \geq m^*$, the variance $\lambda_j(t)$ is bounded from below as follows

$$
\begin{aligned}
\lambda_j^2(t) &\geq \frac{1}{n} \inf_{t \in \mathbb{K}} 2^{jd} \left( m^* \sum_l \phi^2(2^j t - l) - 2^{-jd}((2M)^d \|\phi\|_\infty^2 \|f\|_\infty)^2 \right) \\
&\geq \frac{m^2 m^*}{2}\, \frac{2^{jd}}{n}.
\end{aligned}
$$

Finally, we get

$$E(|Z_{i,n}(t)|^3) \leq c \left( \frac{2^{jd}}{n} \right)^{1/2} \tag{5.11}$$

for $c$ is depending on $m, m^*, \phi$. Let us recall Esseen's inequality

**Proposition 5.1** (Theorem 2.6, Hall (9))**.** *Let* $\{Z_{i,n}, 1 \leq i \leq n\}$ *be a triangular array of independent variables, centered such that* $\sum_{i=1}^n E(Z_{i,n}^2) = 1$. *If there exists some* $\Delta_n \to 0$ *as* $n \to \infty$ *such that*

$$\sum_{i=1}^n E(|Z_{i,n}|^3) \leq \Delta_n,$$

*then, there exists a positive universal constant* $a$ *such that*

$$\forall x \in R, \quad \left| P \left( \sum_{i=1}^n Z_{i,n} \leq x \right) - F_{\mathcal{N}(0,1)}(x) \right| \leq \frac{b}{1+x^2} \Delta_n.$$

We apply this inequality for $\Delta_n = (2^{jd}n^{-1})^{1/2}$, see (5.11), and for each fixed $t \in \mathbb{K}$. We get, for $j, n$ large enough

$$\forall t, \forall x, \quad |F_{\mathcal{X}_j}(x) - F_{\mathcal{N}(0,1)}(x)| \leq \frac{b}{(1+x^2)} \left( \frac{2^{jd}}{n} \right)^{1/2}$$

implying that

$$|z| \quad \leq \quad \pi b \, K(k,j) \left( \frac{2^{jd}}{n} \right)^{1/2}.$$

### 5.2.2. Proof of Lemma 5.2

We proceed as in the proof of Lemma 5.1. Indeed, we have proved there, that for all $k = 1, \ldots, N$ and $t \in \mathbb{K}$

$$\left| K(k,j)E(\cos(\pi k \hat{f}_j(t))) - \cos(\pi k E_j f(t)) \right| \leq K(k,j) u_n,$$

where $|u_n| \leq b 2^{jd} n^{-1/2}$ tends to 0 when $n \to \infty$. Moreover this entails

$$\left| K(k,j)^2 E^2(\cos(\pi k \hat{f}_j(t))) - \cos^2(\pi k E_j f(t)) \right| \leq K(k,j)^2 u_n. \qquad (5.12)$$

For the first term of the needed variance, we have

$$\left| K(k,j)^2 E(\cos^2(\pi k \hat{f}_j(t))) - \frac{1}{2}K(k,j)^2 - K(k,j)^2 \cos(2\pi k E_j f(t))) \right|$$
$$\leq \quad \frac{1}{2} \left| K(k,j)^2 E(\cos(2\pi k \hat{f}_j(t))) - K(k,j)^2 \cos(2\pi k E_j f(t)) \right| \leq \frac{1}{2}K(k,j)^2 u_n.$$

The proof of this last inequality is similar to the proof of Lemma 5.1. Together with (5.12), we get

$$V\left( K(k,j)\cos(\pi k \hat{f}_j(t)) \right)$$
$$\leq \quad \left| K(k,j)^2 E(\cos^2(\pi k \hat{f}_j(t))) - \frac{1}{2}K(k,j)^2 - \frac{1}{2}K(k,j)^2 \cos(2\pi k E_j f(t))) \right|$$
$$+ \left| K(k,j)^2 E^2(\cos(\pi k \hat{f}_j(t))) - \cos^2(\pi k E_j f(t)) \right|$$
$$+ \left| \frac{1}{2}K(k,j)^2 + \frac{1}{2}K(k,j)^{-2}\cos(2\pi k E_j f(t)) - \left( \frac{1}{2} + \frac{1}{2}\cos(2\pi k E_j f(t)) \right) \right|$$
$$\leq \quad \frac{3}{2}K(k,j)^2 u_n + \frac{1}{2} \left| \left( K(k,j)^2 - \cos(2\pi k E_j f(t)) \right) \left( 1 - K(k,j)^{-2} \right) \right|$$
$$\leq \quad (3u_n/2 + \pi^2 k^2 \lambda_j^2) e^{\pi^2 k^2 \lambda_j^2}.$$

Note that $u_n$ is the dominant term for $k$ smaller than $n^{1/4} 2^{-jd/4}$ and $\pi^2 k^2 \lambda_j^2$ is dominant for $k$ larger than the same value.

**Acknowledgements**

The authors thank T. Cai and M. Low for having introduced them to the problem of excess mass estimation.